\documentclass{amsart}
\usepackage{amssymb}
\usepackage{amsthm}

\theoremstyle{definition} \newtheorem{defin}{Definition}[section]
\newtheorem{ex}{Example}[section]
\theoremstyle{remark} \newtheorem{rem}{Remark}
\newcommand{\ab}{\par\medskip\noindent}
\begin{document}

\title{Short catalog of plane ten-edge trees}
\author{Yury Kochetkov}
\date{}
\email{yuyukochetkov@gmail.com, yukochetkov@hse.ru}
\begin{abstract} A type is the set all pairwise nonisotopic plane binary
trees with the same passport. A type is called decomposable, if it
is a union of several Galois orbits. In this work we present the
list of all passports of plane binary trees with ten edges and the
list of all Galois orbits.
\end{abstract} \maketitle

\section{Introduction}

\subsection{Rotation group} A plane tree is a tree imbedded in the plane. Two
trees are considered the same if there exists an isotopy that maps one tree
into another. A plane tree possesses a binary structure --- a coloring of its
vertices into two colors: black and white (adjacent vertices have different
colors). \ab Two operators acts on edges of a plane binary tree: $\alpha$ ---
the counterclockwise rotation around white vertices, $\beta$ --- the
counterclockwise rotation around black vertices. If some enumeration of edges
is given, then operators $\alpha$ and $\beta$ define permutations $a$ and $b$
in $S_n$, where $n$ is the number of edges. Subgroup $\langle a,b\rangle\subset
S_n$ generated by permutations $a$ and $b$ is the \emph{rotation group} of the
tree \cite{LZ}. The rotation group is defined up to conjugation.

\begin{ex} For the tree
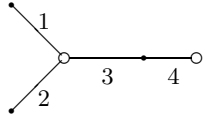
\[\begin{picture}(70,50) \put(0,5){\circle*{2}} \put(0,45){\circle*{2}}
\put(20,25){\circle{4}} \put(50,25){\circle*{2}} \put(70,25){\circle{4}}
\put(0,5){\line(1,1){19}} \put(0,45){\line(1,-1){19}}
\put(22,25){\line(1,0){46}} \put(10,7){\small 2} \put(10,36){\small 1}
\put(34,15){\small 3} \put(59,15){\small 4} \end{picture}\] $a=(1,2,3)$,
$b=(3,4)$ and $\langle a,b\rangle=S(4)$. \end{ex}
\par\noindent
A map of a plane binary tree $T_1$ with set of vertices $V_1$ and set of edges
$E_1$ to a plane binary tree $T_2$ with set of vertices $V_2$ and set of edges
$E_2$ is a map from $V_1$ to $V_2$ such that
\begin{itemize}
    \item white vertices are mapped into white and black --- into black;
    \item adjacent vertices remain adjacent under the mapping (this condition
    defines the map from $E_1$ to $E_2$);
    \item the map commutes with operators $\alpha$ and $\beta$.
\end{itemize} If $E_1|>|E_2|$, then tree $T_1$ is called \emph{reduced} to tree
$T_2$. We will call a tree \emph{non-reducible}, if it cannot be mapped into a
smaller tree.

\begin{ex} Example of a map.
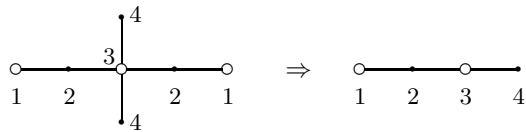
\[\begin{picture}(190,50) \put(0,25){\circle{4}} \put(20,25){\circle*{2}}
\put(40,25){\circle{4}} \put(40,5){\circle*{2}} \put(40,45){\circle*{2}}
\put(60,25){\circle*{2}} \put(80,25){\circle{4}} \put(2,25){\line(1,0){36}}
\put(40,5){\line(0,1){18}} \put(40,45){\line(0,-1){18}}
\put(42,25){\line(1,0){36}} \put(-2,12){\small 1} \put(18,12){\small 2}
\put(33,27){\small 3} \put(43,2){\small 4} \put(43,43){\small 4}
\put(58,12){\small 2} \put(78,12){\small 1} \put(130,25){\circle{4}}
\put(150,25){\circle*{2}} \put(170,25){\circle{4}} \put(190,25){\circle*{2}}
\put(132,25){\line(1,0){36}} \put(172,25){\line(1,0){18}}
\put(102,22){$\Rightarrow$} \put(128,12){\small 1} \put(148,12){\small 2}
\put(168,12){\small 3} \put(188,12){\small 4} \end{picture}\] Here the vertex
with number $i$, $i=1,2,3,4$, is mapped to the vertex with the same number.
Thus, the left six-edge tree is reduced to 3-chain.
\end{ex}

\begin{rem} A tree is reducible, if the set of its edges can be partitioned
into a union of pairwise disjoint classes in a way such that, if two edges
$e_1$ and $e_2$ belong to one class, then edges $\alpha(e_1)$ and $\alpha(e_2)$
belong to one class and edges $\beta(e_1)$ and $\beta(e_2)$ also belong to one
class.
\end{rem}

\begin{rem} Automorphisms of a plane binary tree are rotations around its
center of symmetry (if it exists). Hence, the group of automorphisms is cyclic.
Let $\text{Aut}(T)$ be the group of automorphisms of a tree $T$ and
$|\text{Aut}(T)|$ be the order this group. If the group $\text{Aut}(T)$ is
nontrivial of order $k$, then $T$ will be called $k$-\emph{symmetric}.\end{rem}

\subsection{The Goulden-Jackson formula \cite{GJ}}
\begin{defin} The \emph{passport} of a plane binary tree is the list of
valencies of its white vertices plus the list of valencies of its black
vertices. The set of all trees with the same passport is called a \emph{type}
and is denoted as
$$\Xi=\langle a_1,\ldots,a_m\,|\,b_1,\ldots,b_n\rangle,$$ where $m$ is the
number of white vertices, $n$ is the number of black vertices, $a_1,\ldots,a_m$
is the list of white valencies, written in the nonincreasing order, and
$b_1,\ldots,b_n$ is the list of black valencies (also written in the
nonincreasing order). Thus the tree from Example 1.1 belongs to the type
$\langle 3,1\,|\,2,1,1\rangle$. \end{defin} \ab \textbf{Theorem (the
Goulden-Jackson formula).} \emph{Let $\Xi$ be a type and $w(\Xi)$ be a weighted
sum $w(\Xi)=\sum_{T\in\Xi}(1/|\text{Aut}(T)|)$. Then
$$w(\Xi)=\frac{(m-1)!(n-1)!}{\prod_i k_i!\,\, \prod_j l_j!}.$$ Here $m$ is the number
of white vertices, $n$ is the number of black vertices, $k_i$, $i=1,2,\ldots$,
is the number of white vertices of valency $i$ and $l_j$, $j=1,2,\ldots$, is
the number of black vertices of valency $j$.}

\begin{ex} The Goulden-Jackson formula for the type $\langle
4,1,1\,|\,2,2,1,1\rangle$ gives the number $3/2$. Indeed, two trees belong to
this type
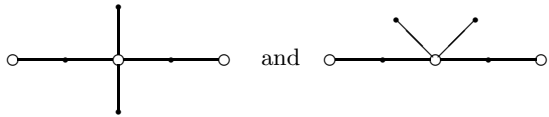
\[\begin{picture}(200,50) \put(0,25){\circle{4}} \put(20,25){\circle*{2}}
\put(40,25){\circle{4}} \put(40,5){\circle*{2}} \put(40,45){\circle*{2}}
\put(60,25){\circle*{2}} \put(80,25){\circle{4}} \put(2,25){\line(1,0){36}}
\put(40,5){\line(0,1){18}} \put(40,27){\line(0,1){18}}
\put(42,25){\line(1,0){36}} \put(94,23){\small and}

\put(120,25){\circle{4}} \put(140,25){\circle*{2}} \put(160,25){\circle{4}}
\put(145,40){\circle*{2}} \put(175,40){\circle*{2}} \put(180,25){\circle*{2}}
\put(200,25){\circle{4}} \put(122,25){\line(1,0){36}}
\put(145,40){\line(1,-1){14}} \put(175,40){\line(-1,-1){14}}
\put(162,25){\line(1,0){36}} \end{picture}\] The left tree is 2-symmetric and
the right is non-symmetric. \end{ex}

\subsection{Galois orbits}

\begin{defin} A polynomial $p\in \mathbb{C}[z]$ is called a \emph{generalized
Chebyshev polynomial} or \emph{Shabat polynomial}, if it has exactly two finite
critical values: 0 and 1. The set $p^{-1}[0,1]$ is a plane binary tree: white
vertices are inverse images of $\{0\}$ and black vertices of $\{1\}$.
\end{defin}
\par\noindent
For each plane binary tree $T$ there exists a Shabat polynomial $p$ with
algebraic coefficients such, that the tree $p^{-1}[0,1]$ is isotopic to $T$. In
what follows all our Shabat polynomials will have algebraic coefficients. \ab
Polynomial $p$ is unique up to linear change of variable $z$. Let $S(t)$ be the
set of all Shabat polynomials that correspond to a tree $T$, i.e. $p\in S(T)
\Rightarrow p^{-1}[0,1]=T$. Each $p\in S(T)$ is defined over some algebraic
field and let $K$ be the minimal of such fields. $K$ will be called the
\emph{definition field} of the tree $T$. If the definition field is
$\mathbb{Q}$, then the tree will be called \emph{rational}. \ab Let $p$ be a
Shabat polynomial with algebraic coefficients and $\gamma$ be an elements of
absolute Galois group. The action of $\gamma$ on $p$ is the action on its
coefficients:
\begin{multline*}q=\gamma(p)=\gamma(a_nz^n+a_{n-1}z^{n-1}+\ldots+a_1z+a_0)=\\=
\gamma(a_n)z^n+\gamma(a_{n-1})z^{n-1}+\ldots+\gamma(a_1)z+\gamma(a_0).
\end{multline*} Then $q$ is a Shabat polynomial and trees $p^{-1}[0,1]$ and
$q^{-1}[0,1]$ belong to the same type. \ab The action of Galois group can be
defined on trees. Indeed, if Shabat polynomials $p_1$ and $p_2$ define equal
trees , then $p_2(z)=p_1(c\,z+d)$, where $c$ and $d$ are algebraic numbers.
Then
$$\gamma(p_2(z))=\gamma(p_1)(\gamma(c)\,z+\gamma(d)),$$ i.e. Shabat polynomials
$\gamma(p_1)$ and $\gamma(p_2)$ also define equal trees. Thus, a type is a
union of disjoint Galois orbits and the degree of definition field of a given
tree is equal to the cardinality of its Galois orbit.

\begin{defin} A type $\Xi$ is called \emph{non-decomposable}, if all its trees
are in one Galois orbit. Otherwise, a type is called
\emph{decomposable}.\end{defin}
\par\noindent
The order of rotation group is a Galois invariant \cite{LZ}. Thus, if a type
contains two trees with different rotation groups, then it is decomposable.
Analogously, a type is decomposable if it contains symmetric trees and
non-symmetric trees.Also a type is decomposable, if it contains two trees such,
that Shabat polynomial of the first is a degree of Shabat polynomial and Shabat
polynomial of the second is not. \ab A decomposable type is called
\emph{trivially decomposable}, if its partition into orbits can be explained by
above-mentioned combinatorial reasons. Otherwise, a type is called
\emph{non-trivially decomposable}. Examples of non-trivial decomposition are
known only for trees of diameter four.

\begin{ex} One of the first discovered example of non-trivial decomposition is
the "Leila flower" \cite{Sc,Ko1}:
$$\Xi=\langle 5,\,\underbrace{1,\ldots,1}_{15}\,|\,6,5,4,3,2\rangle.$$ A tree
in this type is completely defined by cyclic order of black vertices under the
counterclockwise going around of white vertex of valency 5. If the circuit
begins with the vertex of valency 6, then a tree is defined by permutation of
numbers $2,3,4,5$. Here even permutations correspond to trees in one orbit and
odd --- in another. \end{ex}

\section{Setting of the problem}
\[\begin{tabular}{|c|c|c||c|c|c|}\hline N&passport&$w(\Xi)$&N&passport&$w(\Xi)$\\ \hline
1&$\langle 10\,|\,1, 1, 1, 1, 1, 1, 1, 1, 1, 1\rangle$& 1/10 & 2&$\langle 9,
1\,|\,2, 1, 1, 1, 1, 1, 1, 1, 1\rangle$& 1 \\ 3&$\langle 8, 2\,|\,2, 1,
1, 1, 1, 1, 1, 1, 1\rangle$&1 &4&$\langle 8, 1, 1\,|\,3, 1, 1, 1, 1, 1, 1, 1\rangle$&1\\
5&$\langle 8, 1, 1\,|\,2, 2, 1, 1, 1, 1, 1, 1\rangle$&7/2&
6&$\langle 7, 3\,|\,2, 1, 1, 1, 1, 1, 1, 1, 1\rangle$&1\\
7&$\langle 7, 2, 1\,|\,3, 1, 1, 1, 1, 1, 1, 1\rangle$&2&
8&$\langle 7, 2, 1\,|\,2, 2, 1, 1, 1, 1, 1, 1\rangle$&7\\
9&$\langle 7, 1, 1, 1\,|\,4, 1, 1, 1, 1, 1, 1\rangle$&1&
10&$\langle 7, 1, 1, 1\,|\,3, 2, 1, 1, 1, 1, 1\rangle$&6\\
11&$\langle 7, 1, 1, 1\,|\,2, 2, 2, 1, 1, 1, 1\rangle$&5&
12&$\langle 6, 4\,|\,2, 1, 1, 1, 1, 1, 1, 1, 1\rangle$&1\\
13&$\langle 6, 3, 1\,|\,3, 1, 1, 1, 1, 1, 1, 1\rangle$&2&
14&$\langle 6, 3, 1\,|\,2, 2, 1, 1, 1, 1, 1, 1\rangle$&7\\
15&$\langle 6, 2, 2\,|\,3, 1, 1, 1, 1, 1, 1, 1\rangle$&1&
16&$\langle 6, 2, 2\,|\,2, 2, 1, 1, 1, 1, 1, 1\rangle$&7/2\\
17&$\langle 6, 2, 1, 1\,|\,4, 1, 1, 1, 1, 1, 1\rangle$&3&
18&$\langle 6, 2, 1, 1\,|\,3, 2, 1, 1, 1, 1, 1\rangle$&18\\
19&$\langle 6, 2, 1, 1\,|\,2, 2, 2, 1, 1, 1, 1\rangle$&15&
20&$\langle 6, 1, 1, 1, 1\,|\,5, 1, 1, 1, 1, 1\rangle$&1\\
21&$\langle 6, 1, 1, 1, 1\,|\,4, 2, 1, 1, 1, 1\rangle$&5&
22&$\langle 6, 1, 1, 1, 1\,|\,3, 3, 1, 1, 1, 1\rangle$&5/2\\
23&$\langle 6, 1, 1, 1, 1\,|\,3, 2, 2, 1, 1, 1\rangle$&10&
24&$\langle 6, 1, 1, 1, 1\,|\,2, 2, 2, 2, 1, 1\rangle$&5/2\\
25&$\langle 5, 5\,|\,2, 1, 1, 1, 1, 1, 1, 1, 1\rangle$&1/2&
26&$\langle 5, 4, 1\,|\,3, 1, 1, 1, 1, 1, 1, 1\rangle$&2\\
27&$\langle 5, 4, 1\,|\,2, 2, 1, 1, 1, 1, 1, 1\rangle$&7&
28&$\langle 5, 3, 2\,|\,3, 1, 1, 1, 1, 1, 1, 1\rangle$&2\\
29&$\langle 5, 3, 2\,|\,2, 2, 1, 1, 1, 1, 1, 1\rangle$&7&
30&$\langle 5, 3, 1, 1\,|\,4, 1, 1, 1, 1, 1, 1\rangle$&3\\
31&$\langle 5, 3, 1, 1\,|\,3, 2, 1, 1, 1, 1, 1\rangle$&18&
32&$\langle 5, 3, 1, 1\,|\,2, 2, 2, 1, 1, 1, 1\rangle$&15\\
33&$\langle 5, 2, 2, 1\,|\,4, 1, 1, 1, 1, 1, 1\rangle$&3&
34&$\langle 5, 2, 2, 1\,|\,3, 2, 1, 1, 1, 1, 1\rangle$&18\\
35&$\langle 5, 2, 2, 1\,|\,2, 2, 2, 1, 1, 1, 1\rangle$&15&
36&$\langle 5, 2, 1, 1, 1\,|\,5, 1, 1, 1, 1, 1\rangle$&4\\
37&$\langle 5, 2, 1, 1, 1\,|\,4, 2, 1, 1, 1, 1\rangle$&20&
38&$\langle 5, 2, 1, 1, 1\,|\,3, 3, 1, 1, 1, 1\rangle$&10\\
39&$\langle 5, 2, 1, 1, 1\,|\,3, 2, 2, 1, 1, 1\rangle$&40&
40&$\langle 5, 2, 1, 1, 1\,|\,2, 2, 2, 2, 1, 1\rangle$&10\\
41&$\langle 5, 1, 1, 1, 1, 1\,|\,4, 3, 1, 1, 1\rangle$&4&
42&$\langle 5, 1, 1, 1, 1, 1\,|\,4, 2, 2, 1, 1\rangle$&6\\
43&$\langle 5, 1, 1, 1, 1, 1\,|\,3, 3, 2, 1, 1\rangle$&6&
44&$\langle 5, 1, 1, 1, 1, 1\,|\,3, 2, 2, 2, 1\rangle$&4\\
45&$\langle 5, 1, 1, 1, 1, 1\,|\,2, 2, 2, 2, 2\rangle$&1/5&
46&$\langle 4, 4, 2\,|\,3, 1, 1, 1, 1, 1, 1, 1\rangle$&1\\
47&$\langle 4, 4, 2\,|\,2, 2, 1, 1, 1, 1, 1, 1\rangle$&7/2&
48&$\langle 4, 4, 1, 1\,|\,4, 1, 1, 1, 1, 1, 1\rangle$&3/2\\
49&$\langle 4, 4, 1, 1\,|\,3, 2, 1, 1, 1, 1, 1\rangle$&9&
50&$\langle 4, 4, 1, 1\,|\,2, 2, 2, 1, 1, 1, 1\rangle$&15/2\\
51&$\langle 4, 3, 3\,|\,3, 1, 1, 1, 1, 1, 1, 1\rangle$&1&
52&$\langle 4, 3, 3\,|\,2, 2, 1, 1, 1, 1, 1, 1\rangle$&7/2\\
53&$\langle 4, 3, 2, 1\,|\,4, 1, 1, 1, 1, 1, 1\rangle$&6&
54&$\langle 4, 3, 2, 1\,|\,3, 2, 1, 1, 1, 1, 1\rangle$&36\\
55&$\langle 4, 3, 2, 1\,|\,2, 2, 2, 1, 1, 1, 1\rangle$&30&
56&$\langle 4, 3, 1, 1, 1\,|\,4, 2, 1, 1, 1, 1\rangle$&20\\
57&$\langle 4, 3, 1, 1, 1\,|\,3, 3, 1, 1, 1, 1\rangle$&10&
58&$\langle 4, 3, 1, 1, 1\,|\,3, 2, 2, 1, 1, 1\rangle$&40\\
59&$\langle 4, 3, 1, 1, 1\,|\,2, 2, 2, 2, 1, 1\rangle$&10&
60&$\langle 4, 2, 2, 2\,|\,4, 1, 1, 1, 1, 1, 1\rangle$&1\\
61&$\langle 4, 2, 2, 2\,|\,3, 2, 1, 1, 1, 1, 1\rangle$&6&
62&$\langle 4, 2, 2, 2\,|\,2, 2, 2, 1, 1, 1, 1\rangle$&5\\
63&$\langle 4, 2, 2, 1, 1\,|\,4, 2, 1, 1, 1, 1\rangle$&30&
64&$\langle 4, 2, 2, 1, 1\,|\,3, 3, 1, 1, 1, 1\rangle$&15\\
65&$\langle 4, 2, 2, 1, 1\,|\,3, 2, 2, 1, 1, 1\rangle$&60&
66&$\langle 4, 2, 2, 1, 1\,|\,2, 2, 2, 2, 1, 1\rangle$&15\\
67&$\langle 4, 2, 1, 1, 1, 1\,|\,3, 3, 2, 1, 1\rangle$&30&
68&$\langle 4, 2, 1, 1, 1, 1\,|\,3, 2, 2, 2, 1\rangle$&20\\
69&$\langle 4, 2, 1, 1, 1, 1\,|\,2, 2, 2, 2, 2\rangle$&1&
70&$\langle 4, 1, 1, 1, 1, 1, 1\,|\,3, 3, 3, 1\rangle$&1\\
71&$\langle 4, 1, 1, 1, 1, 1, 1\,|\,3, 3, 2, 2\rangle$&3/2&
72&$\langle 3, 3, 3, 1\,|\,3, 2, 1, 1, 1, 1, 1\rangle$&6\\
73&$\langle 3, 3, 3, 1\,|\,2, 2, 2, 1, 1, 1, 1\rangle$&5&
74&$\langle 3, 3, 2, 2\,|\,3, 2, 1, 1, 1, 1, 1\rangle$&9\\
75&$\langle 3, 3, 2, 2\,|\,2, 2, 2, 1, 1, 1, 1\rangle$&15/2&
76&$\langle 3, 3, 2, 1, 1\,|\,3, 3, 1, 1, 1, 1\rangle$&15\\
77&$\langle 3, 3, 2, 1, 1\,|\,3, 2, 2, 1, 1, 1\rangle$&60&
78&$\langle 3, 3, 2, 1, 1\,|\,2, 2, 2, 2, 1, 1\rangle$&15\\
79&$\langle 3, 3, 1, 1, 1, 1\,|\,3, 2, 2, 2, 1\rangle$&10&
80&$\langle 3, 3, 1,1, 1, 1\,|\,2, 2, 2, 2, 2\rangle$&1/2\\
81&$\langle 3, 2, 2, 2, 1\,|\,3, 2, 2, 1, 1, 1\rangle$&40&
82&$\langle 3, 2, 2, 2, 1\,|\,2, 2, 2, 2, 1, 1\rangle$&10\\
83&$\langle 3, 2, 2, 1, 1, 1\,|\,2, 2, 2, 2, 2\rangle$&2&
84&$\langle 2, 2, 2, 2, 2\,|\,2, 2, 2, 2, 1, 1\rangle$&1/2\\
\hline \end{tabular}\] \ab In table above all types of ten-edge trees are
enumerated. The lists of white and black valencies are partitions of number 10.
The list of white valencies in each type is lexicographically higher, than the
list of black (as the number of vertices is 11, then these two lists cannot be
equal). For each type we compute the corresponding weighted sum with the use of
Goulden-Jackson formula. This weighted sum is equal to the type cardinality, if
there is no symmetric trees in the type. Otherwise, it is less.

\begin{rem} In what follows we will write "Galois $n$-orbit" instead of
"Galois orbit of cardinality $n$".\end{rem}

\begin{rem}  Shabat polynomials and definition fields of nine-edge trees are
described in catalog \cite{Ko2} and with number of edges $\leqslant 8$ --- in
catalog \cite{Zv}.\end{rem}

\begin{rem} In what follows critical values of Shabat polynomial will be $0$
and some nonzero rational number. \end{rem}

\begin{rem} In what follows we will give only schematic pictures of trees and
will not try to present their true forms (see \cite{Zv, Ko3}).\end{rem}

\section{orbits}

\subsection{Types 1, 25, 45, 80 and 84} Each of them  contains only one tree ---
a symmetric tree with symmetry of the order 10, 2, 5, 2 and 2, respectively.

\subsection{Types 2 --- 4, 6 --- 15, 17 --- 21, 23, 26 ---
32, 34 --- 42, 44, 46, 49, 51, 53 --- 60, 62, 63, 65, 67 --- 70, 72 --- 74, 77,
79, 81 and 82}  All these types are non-decomposable and do not contain
symmetric trees.

\subsection{Types 5, 22, 24, 48, 52 and 71 } Each of them contains one
2-symmetric tree (Galois 1-orbit). Non-symmetric trees in each of these types
constitute Galois 3-orbit, 2-orbit, 2-orbit, 1-orbit, 3-orbit and 1-orbit,
respectively. Hence, non-symmetric trees in types 48 and 71 are rational.
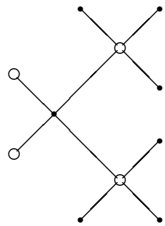
\[\begin{picture}(235,90) \put(0,30){\circle{4}}
\put(0,60){\circle{4}} \put(15,45){\circle*{2}} \put(40,70){\circle{4}}
\put(40,20){\circle{4}} \put(25,5){\circle*{2}} \put(25,85){\circle*{2}}
\put(55,5){\circle*{2}} \put(55,85){\circle*{2}} \put(55,35){\circle*{2}}
\put(55,55){\circle*{2}} \put(1,31){\line(1,1){38}} \put(1,59){\line(1,-1){38}}
\put(39,19){\line(-1,-1){14}} \put(41,19){\line(1,-1){14}}
\put(41,21){\line(1,1){14}} \put(39,71){\line(-1,1){14}}
\put(41,69){\line(1,-1){14}} \put(41,71){\line(1,1){14}}
\put(110,42){$\left(x^2+\frac{100}{27}\right)^4\left(x^2+
\frac{100}{27}\,x+\frac{100}{27}\right)$} \end{picture}\]
\begin{center}{The rational non-symmetric tree in the type 48
and its Shabat polynomial.}\end{center}
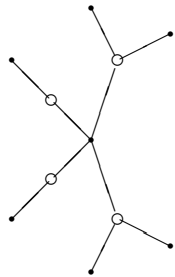
\[\begin{picture}(240,110) \put(0,85){\circle*{2}} \put(0,25){\circle*{2}}
\put(15,70){\circle{4}}\put(15,40){\circle{4}} \put(30,55){\circle*{2}}
\put(40,85){\circle{4}} \put(40,25){\circle{4}} \put(30,105){\circle*{2}}
\put(30,5){\circle*{2}} \put(60,95){\circle*{2}} \put(60,15){\circle*{2}}
\put(0,85){\line(1,-1){14}} \put(0,25){\line(1,1){14}}
\put(30,55){\line(-1,1){14}} \put(30,55){\line(-1,-1){14}}
\put(30,55){\line(1,3){9}} \put(30,55){\line(1,-3){9}}
\put(30,105){\line(1,-2){9}} \put(30,5){\line(1,2){9}}
\put(60,95){\line(-2,-1){19}} \put(60,15){\line(-2,1){19}}
\put(100,52){$\left(x^2+\frac{6075}{5476}\right)^3\left(x^2+\frac{2025}{592}\,x+
\frac{2025}{592}\right)^2$}\end{picture}\]
\begin{center}{The rational non-symmetric tree in the type 71
and its Shabat polynomial.}\end{center}

\subsection{Types 64, 66, 76 and 78} Each of them contains
two 2-symmetric trees and 14 non-symmetric. Symmetric trees constitute Galois
2-orbit and non-symmetric --- 14-orbit.
\[\begin{picture}(210,55) \put(0,35){\circle*{2}} \put(15,35){\circle{4}}
\put(30,35){\circle*{2}} \put(30,20){\circle{4}} \put(45,35){\circle{4}}
\put(45,20){\circle*{2}} \put(45,50){\circle*{2}} \put(60,35){\circle*{2}}
\put(60,50){\circle{4}} \put(75,35){\circle{4}} \put(90,35){\circle*{2}}
\put(0,35){\line(1,0){13}} \put(17,35){\line(1,0){26}}
\put(30,22){\line(0,1){13}} \put(45,20){\line(0,1){13}}
\put(45,37){\line(0,1){13}} \put(47,35){\line(1,0){26}}
\put(60,35){\line(0,1){13}} \put(77,35){\line(1,0){13}}

\put(120,35){\circle*{2}} \put(135,35){\circle{4}} \put(150,35){\circle*{2}}
\put(150,50){\circle{4}} \put(165,35){\circle{4}} \put(165,20){\circle*{2}}
\put(165,50){\circle*{2}} \put(180,35){\circle*{2}} \put(180,20){\circle{4}}
\put(195,35){\circle{4}} \put(210,35){\circle*{2}} \put(120,35){\line(1,0){13}}
\put(137,35){\line(1,0){26}} \put(150,35){\line(0,1){13}}
\put(165,20){\line(0,1){13}} \put(165,37){\line(0,1){13}}
\put(167,35){\line(1,0){26}} \put(180,22){\line(0,1){13}}
\put(197,35){\line(1,0){13}}

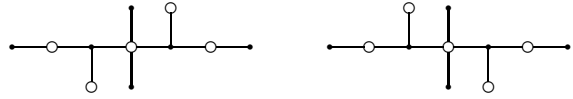
\put(45,2){\small Symmetric trees in the type 64} \end{picture}\]

\[\begin{picture}(240,70) \put(0,35){\circle{4}} \put(15,35){\circle*{2}}
\put(30,35){\circle{4}} \put(45,35){\circle*{2}} \put(60,35){\circle{4}}
\put(60,20){\circle*{2}} \put(60,50){\circle*{2}} \put(75,35){\circle*{2}}
\put(90,35){\circle{4}} \put(105,35){\circle*{2}} \put(120,35){\circle{4}}
\put(2,35){\line(1,0){26}} \put(32,35){\line(1,0){26}}
\put(62,35){\line(1,0){26}} \put(92,35){\line(1,0){26}}
\put(60,20){\line(0,1){13}} \put(60,37){\line(0,1){13}}

\put(150,35){\circle*{2}} \put(165,35){\circle{4}} \put(180,35){\circle*{2}}
\put(195,35){\circle{4}} \put(195,20){\circle*{2}} \put(195,50){\circle*{2}}
\put(195,5){\circle{4}} \put(195,65){\circle{4}} \put(210,35){\circle*{2}}
\put(225,35){\circle{4}} \put(240,35){\circle*{2}} \put(150,35){\line(1,0){13}}
\put(167,35){\line(1,0){26}} \put(195,7){\line(0,1){26}}
\put(195,37){\line(0,1){26}} \put(197,35){\line(1,0){26}}
\put(227,35){\line(1,0){13}}

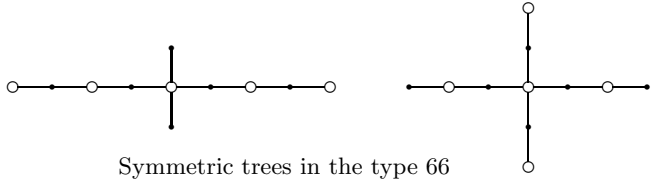
\put(40,2){\small Symmetric trees in the type 66} \end{picture}\]

\[\begin{picture}(210,55) \put(0,35){\circle*{2}} \put(15,35){\circle{4}}
\put(15,50){\circle*{2}} \put(30,35){\circle*{2}} \put(30,20){\circle{4}}
\put(45,35){\circle{4}} \put(60,35){\circle*{2}} \put(60,50){\circle{4}}
\put(75,35){\circle{4}} \put(75,20){\circle*{2}} \put(90,35){\circle*{2}}
\put(0,35){\line(1,0){13}} \put(15,37){\line(0,1){13}}
\put(17,35){\line(1,0){26}} \put(30,35){\line(0,-1){13}}
\put(47,35){\line(1,0){26}} \put(60,35){\line(0,1){13}}
\put(75,33){\line(0,-1){13}} \put(77,35){\line(1,0){13}}

\put(120,35){\circle*{2}} \put(135,35){\circle{4}} \put(135,20){\circle*{2}}
\put(150,35){\circle*{2}} \put(150,50){\circle{4}} \put(165,35){\circle{4}}
\put(180,35){\circle*{2}} \put(180,20){\circle{4}} \put(195,35){\circle{4}}
\put(195,50){\circle*{2}} \put(210,35){\circle*{2}}
\put(120,35){\line(1,0){13}} \put(135,33){\line(0,-1){13}}
\put(137,35){\line(1,0){26}} \put(150,35){\line(0,1){13}}
\put(167,35){\line(1,0){26}} \put(180,35){\line(0,-1){13}}
\put(195,37){\line(0,1){13}} \put(197,35){\line(1,0){13}} \put(45,2){\small
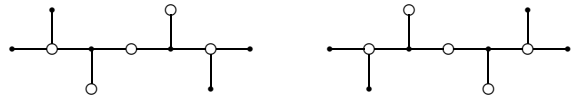
Symmetric trees in the type 76} \end{picture}\]

\[\begin{picture}(270,55) \put(0,35){\circle{4}} \put(15,35){\circle*{2}}
\put(30,35){\circle{4}} \put(30,50){\circle*{2}} \put(45,35){\circle*{2}}
\put(60,35){\circle{4}} \put(75,35){\circle*{2}} \put(90,35){\circle{4}}
\put(90,20){\circle*{2}} \put(105,35){\circle*{2}} \put(120,35){\circle{4}}
\put(2,35){\line(1,0){26}} \put(32,35){\line(1,0){26}}
\put(30,37){\line(0,1){13}} \put(62,35){\line(1,0){26}}
\put(92,35){\line(1,0){26}} \put(90,33){\line(0,-1){13}}

\put(150,35){\circle{4}} \put(165,35){\circle*{2}} \put(180,35){\circle{4}}
\put(180,20){\circle*{2}} \put(195,35){\circle*{2}} \put(210,35){\circle{4}}
\put(225,35){\circle*{2}} \put(240,35){\circle{4}} \put(240,50){\circle*{2}}
\put(255,35){\circle*{2}} \put(270,35){\circle{4}} \put(152,35){\line(1,0){26}}
\put(182,35){\line(1,0){26}} \put(180,33){\line(0,-1){13}}
\put(212,35){\line(1,0){26}} \put(242,35){\line(1,0){26}}
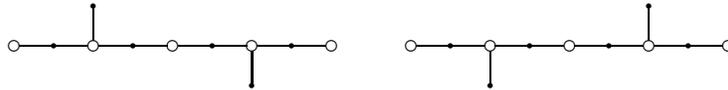
\put(240,37){\line(0,1){13}} \put(75,2){\small Symmetric trees in the type 78}
\end{picture}\]

\subsection{Type 75} It contains three 2-symmetric trees (Galois
3-orbit) and six non-symmetric (Galois 6-orbit). Schemas of symmetric trees are
presented below.
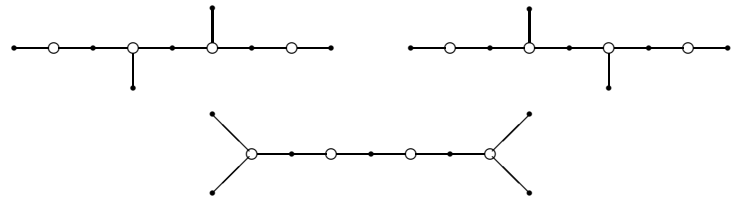
\[\begin{picture}(270,80) \put(0,60){\circle*{2}} \put(15,60){\circle{4}}
\put(30,60){\circle*{2}} \put(45,60){\circle{4}} \put(45,45){\circle*{2}}
\put(60,60){\circle*{2}} \put(75,60){\circle{4}} \put(75,75){\circle*{2}}
\put(90,60){\circle*{2}} \put(105,60){\circle{4}} \put(120,60){\circle*{2}}
\put(0,60){\line(1,0){13}} \put(17,60){\line(1,0){26}}
\put(47,60){\line(1,0){26}} \put(77,60){\line(1,0){26}}
\put(107,60){\line(1,0){13}} \put(45,45){\line(0,1){13}}
\put(75,75){\line(0,-1){13}}

\put(150,60){\circle*{2}} \put(165,60){\circle{4}} \put(180,60){\circle*{2}}
\put(195,60){\circle{4}} \put(195,75){\circle*{2}} \put(210,60){\circle*{2}}
\put(225,60){\circle{4}} \put(225,45){\circle*{2}} \put(240,60){\circle*{2}}
\put(255,60){\circle{4}} \put(270,60){\circle*{2}} \put(150,60){\line(1,0){13}}
\put(167,60){\line(1,0){26}} \put(197,60){\line(1,0){26}}
\put(227,60){\line(1,0){26}} \put(257,60){\line(1,0){13}}
\put(195,75){\line(0,-1){13}} \put(225,45){\line(0,1){13}}

\put(75,35){\circle*{2}} \put(75,5){\circle*{2}} \put(90,20){\circle{4}}
\put(105,20){\circle*{2}} \put(120,20){\circle{4}} \put(135,20){\circle*{2}}
\put(150,20){\circle{4}} \put(165,20){\circle*{2}} \put(180,20){\circle{4}}
\put(195,5){\circle*{2}} \put(195,35){\circle*{2}} \put(75,35){\line(1,-1){14}}
\put(75,5){\line(1,1){14}} \put(92,20){\line(1,0){26}}
\put(122,20){\line(1,0){26}} \put(152,20){\line(1,0){26}}
\put(195,35){\line(-1,-1){14}} \put(195,5){\line(-1,1){14}} \end{picture}\]

\subsection{Type 16}  It contains four trees. Two of them --- $T_1$
and $T_2$
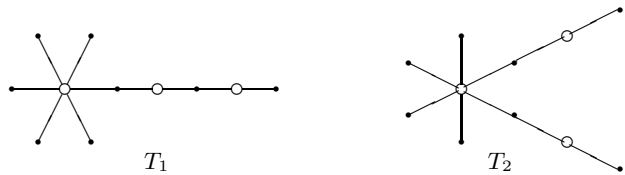
\[\begin{picture}(230,70) \put(0,35){\circle*{2}} \put(10,55){\circle*{2}}
\put(10,15){\circle*{2}} \put(20,35){\circle{4}} \put(30,55){\circle*{2}}
\put(30,15){\circle*{2}} \put(40,35){\circle*{2}} \put(55,35){\circle{4}}
\put(70,35){\circle*{2}} \put(85,35){\circle{4}} \put(100,35){\circle*{2}}
\put(0,35){\line(1,0){18}} \put(10,55){\line(1,-2){9}}
\put(10,15){\line(1,2){9}} \put(30,55){\line(-1,-2){9}}
\put(30,15){\line(-1,2){9}} \put(22,35){\line(1,0){31}}
\put(57,35){\line(1,0){26}} \put(87,35){\line(1,0){13}} \put(50,5){\small
$T_1$}

\put(150,45){\circle*{2}} \put(150,25){\circle*{2}} \put(170,35){\circle{4}}
\put(170,55){\circle*{2}} \put(170,15){\circle*{2}} \put(190,45){\circle*{2}}
\put(190,25){\circle*{2}} \put(210,55){\circle{4}} \put(210,15){\circle{4}}
\put(230,65){\circle*{2}} \put(230,5){\circle*{2}}
\put(150,45){\line(2,-1){19}} \put(150,25){\line(2,1){19}}
\put(170,15){\line(0,1){18}} \put(170,55){\line(0,-1){18}}
\put(171,36){\line(2,1){38}} \put(171,34){\line(2,-1){38}}
\put(211,56){\line(2,1){18}} \put(211,14){\line(2,-1){18}} \put(180,5){\small
$T_2$} \end{picture}\] constitute 2-orbit. They can be reduced to 2-chain and
their rotation group has order 14400.\ab 2-symmetric tree $T_3$
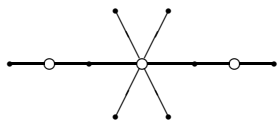
\[\begin{picture}(100,50) \put(0,25){\circle*{2}} \put(15,25){\circle{4}}
\put(30,25){\circle*{2}} \put(40,45){\circle*{2}} \put(40,5){\circle*{2}}
\put(50,25){\circle{4}} \put(60,45){\circle*{2}} \put(60,5){\circle*{2}}
\put(70,25){\circle*{2}} \put(85,25){\circle{4}} \put(100,25){\circle*{2}}
\put(0,25){\line(1,0){13}} \put(17,25){\line(1,0){31}}
\put(40,45){\line(1,-2){9}} \put(40,5){\line(1,2){9}}
\put(60,45){\line(-1,-2){9}} \put(60,5){\line(-1,2){9}}
\put(52,25){\line(1,0){31}} \put(87,25){\line(1,0){13}} \end{picture}\]
constitute 1-orbit. The order of its rotation group is 240. \ab Tree $T_4$
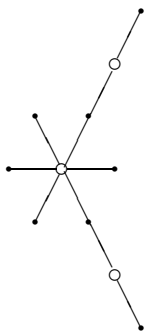
\[\begin{picture}(50,130) \put(0,65){\circle*{2}}
\put(10,85){\circle*{2}} \put(10,45){\circle*{2}} \put(20,65){\circle{4}}
\put(30,85){\circle*{2}} \put(30,45){\circle*{2}} \put(40,65){\circle*{2}}
\put(40,105){\circle{4}} \put(40,25){\circle{4}} \put(50,125){\circle*{2}}
\put(50,5){\circle*{2}} \put(0,65){\line(1,0){18}} \put(10,85){\line(1,-2){9}}
\put(10,45){\line(1,2){9}} \put(22,65){\line(1,0){18}}
\put(21,66){\line(1,2){18}} \put(21,64){\line(1,-2){18}}
\put(50,125){\line(-1,-2){9}} \put(50,5){\line(-1,2){9}} \end{picture}\] is
rational and constitute 1-orbit. Its Shabat polynomial $x^6(x^2-2x+32/5)^2$ is
the square of Shabat polynomial $x^3(x^2-2x+32/5)$ that corresponds to the tree
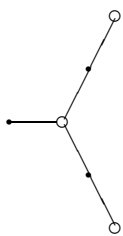
\[\begin{picture}(40,95) \put(0,45){\circle*{2}}
\put(20,45){\circle{4}} \put(30,25){\circle*{2}} \put(30,65){\circle*{2}}
\put(40,5){\circle{4}} \put(40,85){\circle{4}} \put(0,45){\line(1,0){18}}
\put(21,46){\line(1,2){19}} \put(21,44){\line(1,-2){19}}\end{picture}\] In
"squaring a tree", inverse images of segment $[-1,0]$ become new edges. \ab The
tree $T_4$ can be reduced to 2-chain. The order of its rotation group is 7200.

\subsection{Type 47} Here everything is more or less the same, as in the
previous type 16. There are also 4 trees. Two of them belong to the Galois
2-orbit. They can be reduced to 2-chain and their rotation group has order
14400. One tree is 2-symmetric
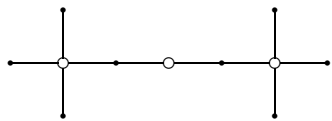
\[\begin{picture}(120,50) \put(0,25){\circle*{2}} \put(20,25){\circle{4}}
\put(20,5){\circle*{2}} \put(20,45){\circle*{2}} \put(40,25){\circle*{2}}
\put(60,25){\circle{4}} \put(80,25){\circle*{2}} \put(100,25){\circle{4}}
\put(100,5){\circle*{2}} \put(100,45){\circle*{2}} \put(120,25){\circle*{2}}
\put(0,25){\line(1,0){18}} \put(20,5){\line(0,1){18}}
\put(20,45){\line(0,-1){18}} \put(22,25){\line(1,0){36}}
\put(62,25){\line(1,0){36}} \put(100,5){\line(0,1){18}}
\put(100,45){\line(0,-1){18}} \put(102,25){\line(1,0){18}} \end{picture}\] Its
rotation group has order 240. \ab The fourth tree
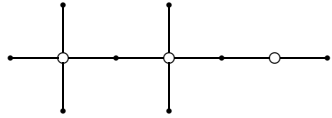
\[\begin{picture}(120,60) \put(0,25){\circle*{2}} \put(20,25){\circle{4}}
\put(20,5){\circle*{2}} \put(20,45){\circle*{2}} \put(40,25){\circle*{2}}
\put(60,25){\circle{4}} \put(60,5){\circle*{2}} \put(60,45){\circle*{2}}
\put(80,25){\circle*{2}} \put(100,25){\circle{4}} \put(120,25){\circle*{2}}
\put(0,25){\line(1,0){18}} \put(20,5){\line(0,1){18}}
\put(20,45){\line(0,-1){18}} \put(22,25){\line(1,0){36}}
\put(62,25){\line(1,0){36}} \put(60,5){\line(0,1){18}}
\put(60,45){\line(0,-1){18}} \put(102,25){\line(1,0){18}}\end{picture}\] is
rational and can be reduced to 2-chain. The order of its rotation group is 200.
Its Shabat polynomial $(x^2-1/5)^4(x-1)^2$ is the square of Shabat polynomial
$(x^2-1/5)^2(x-1)$ that corresponds to 5-chain
\[\begin{picture}(100,10) \put(0,5){\circle*{2}} \put(20,5){\circle{4}}
\put(40,5){\circle*{2}} \put(60,5){\circle{4}} \put(80,5){\circle*{2}}
\put(100,5){\circle{4}} \put(0,5){\line(1,0){18}} \put(22,5){\line(1,0){36}}
\put(62,5){\line(1,0){36}}\end{picture}\]

\subsection{Type 61}  It contains 6 trees. Each of them can be reduced to
2-chain. The rotation group of each tree has order 28800. The type is
decomposable: 5 trees constitute Galois 5-orbit and the sixth tree
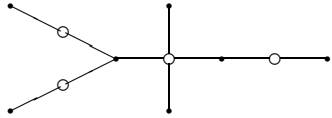
\[\begin{picture}(120,50) \put(0,45){\circle*{2}} \put(0,5){\circle*{2}}
\put(20,35){\circle{4}} \put(20,15){\circle{4}} \put(40,25){\circle*{2}}
\put(60,25){\circle{4}} \put(60,45){\circle*{2}} \put(60,5){\circle*{2}}
\put(80,25){\circle*{2}} \put(100,25){\circle{4}} \put(120,25){\circle*{2}}
\put(0,45){\line(2,-1){19}} \put(0,5){\line(2,1){19}}
\put(40,25){\line(-2,1){19}} \put(40,25){\line(-2,-1){19}}
\put(40,25){\line(1,0){18}} \put(60,5){\line(0,1){18}}
\put(60,45){\line(0,-1){18}} \put(62,25){\line(1,0){36}}
\put(102,25){\line(1,0){18}}\end{picture}\] is rational. Its Shabat polynomial
$$x^4\left(x^3+\frac 59 x^2-\frac{5}{81} x -\frac{5}{81}\right)^2$$ is the
square of Shabat polynomial that corresponds to the tree
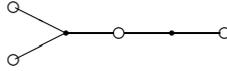
\[\begin{picture}(80,30) \put(0,25){\circle{4}}  \put(0,5){\circle{4}}
\put(20,15){\circle*{2}}  \put(40,15){\circle{4}}  \put(60,15){\circle*{2}}
\put(80,15){\circle{4}} \put(20,15){\line(-2,1){19}}
\put(20,15){\line(-2,-1){19}} \put(20,15){\line(1,0){18}}
\put(42,15){\line(1,0){36}} \end{picture}\]

\subsection{Type 33}  It contains 3 trees. The rotation group of each
of them is $S_{10}$. Two trees belong to Galois 2-orbit and the tree
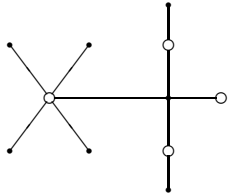
\[\begin{picture}(100,80) \put(0,60){\circle*{2}} \put(0,20){\circle*{2}}
\put(15,40){\circle{4}} \put(0,60){\line(3,-4){14}} \put(0,20){\line(3,4){14}}
\put(30,60){\circle*{2}}  \put(30,20){\circle*{2}}
\put(30,60){\line(-3,-4){14}} \put(30,20){\line(-3,4){14}}
\put(60,40){\circle*{2}} \put(17,40){\line(1,0){43}} \put(60,60){\circle{4}}
\put(60,20){\circle{4}} \put(60,75){\circle*{2}} \put(60,5){\circle*{2}}
\put(80,40){\circle{4}} \put(60,40){\line(1,0){18}} \put(60,40){\line(0,1){18}}
\put(60,40){\line(0,-1){18}} \put(60,5){\line(0,1){13}}
\put(60,75){\line(0,-1){13}} \end{picture}\] is rational with Shabat polynomial
$$x^5\left(x^2-\frac{27}{16}\,+\frac{27}{32}\right)^2(x-1).$$ Trees in this
type are trees of diameter 4. The black vertex of valency 4 is the center and a
cyclic order of white vertices under counterclockwise going around of the
center defines the tree. This type is an example of non-trivial decomposition.

\subsection{Type 43}  It contains 6 trees. The rotation group of each
of them is $S_{10}$. Five trees constitute a Galois 5-orbit and the sixth tree
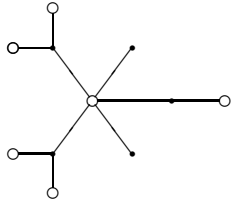
\[\begin{picture}(80,80) \put(0,60){\circle{4}} \put(0,20){\circle{4}}
\put(15,60){\circle*{2}} \put(15,20){\circle*{2}} \put(15,75){\circle{4}}
\put(15,5){\circle{4}} \put(0,60){\circle{4}} \put(30,40){\circle{4}}
\put(45,60){\circle*{2}} \put(45,20){\circle*{2}} \put(60,40){\circle*{2}}
\put(80,40){\circle{4}} \put(2,60){\line(1,0){13}}  \put(2,20){\line(1,0){13}}
\put(15,60){\line(0,1){13}}  \put(15,20){\line(0,-1){13}}
\put(15,60){\line(3,-4){14}}  \put(15,20){\line(3,4){14}}
\put(45,60){\line(-3,-4){14}}  \put(45,20){\line(-3,4){14}}
\put(32,40){\line(1,0){46}} \end{picture}\] is rational with Shabat polynomial
$$\int x^4(x-1)(x^2+2x+2)^2 d\,x.$$  Trees in this
type are trees of diameter 4. The white vertex of valency 5 is the center and a
cyclic order of black vertices under counterclockwise going around of the
center defines the tree. This type also is an example of non-trivial
decomposition.

\subsection{Type 50}  It contains 9 trees. Three of them are
2-symmetric
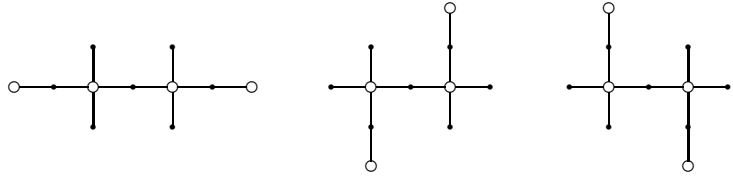
\[\begin{picture}(270,70) \put(0,35){\circle{4}} \put(15,35){\circle*{2}}
\put(30,35){\circle{4}} \put(30,20){\circle*{2}} \put(30,50){\circle*{2}}
\put(45,35){\circle*{2}} \put(60,35){\circle{4}} \put(60,20){\circle*{2}}
\put(60,50){\circle*{2}} \put(75,35){\circle*{2}} \put(90,35){\circle{4}}
\put(2,35){\line(1,0){26}} \put(30,20){\line(0,1){13}}
\put(30,50){\line(0,-1){13}} \put(32,35){\line(1,0){26}}
\put(60,20){\line(0,1){13}} \put(60,50){\line(0,-1){13}}
\put(62,35){\line(1,0){26}}

\put(120,35){\circle*{2}} \put(135,35){\circle{4}} \put(135,20){\circle*{2}}
\put(135,5){\circle{4}} \put(135,50){\circle*{2}} \put(150,35){\circle*{2}}
\put(165,35){\circle{4}} \put(165,50){\circle*{2}} \put(165,65){\circle{4}}
\put(165,20){\circle*{2}} \put(180,35){\circle*{2}}
\put(120,35){\line(1,0){13}} \put(135,33){\line(0,-1){26}}
\put(135,37){\line(0,1){13}} \put(137,35){\line(1,0){26}}
\put(165,37){\line(0,1){26}} \put(165,33){\line(0,-1){13}}
\put(167,35){\line(1,0){13}}

\put(210,35){\circle*{2}} \put(225,35){\circle{4}} \put(225,50){\circle*{2}}
\put(225,65){\circle{4}} \put(225,20){\circle*{2}} \put(240,35){\circle*{2}}
\put(255,35){\circle{4}} \put(255,20){\circle*{2}} \put(255,5){\circle{4}}
\put(255,50){\circle*{2}} \put(270,35){\circle*{2}}
\put(210,35){\line(1,0){13}} \put(225,37){\line(0,1){26}}
\put(225,33){\line(0,-1){13}} \put(227,35){\line(1,0){26}}
\put(255,33){\line(0,-1){26}} \put(255,37){\line(0,1){13}}
\put(257,35){\line(1,0){13}} \end{picture}\] and they constitute Galois
3-orbit. The rotation group of each of them has order 3840. \ab One tree
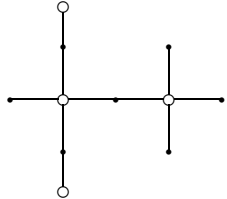
\[\begin{picture}(80,80) \put(0,40){\circle*{2}} \put(20,40){\circle{4}} \put(20,60){\circle*{2}}
\put(20,20){\circle*{2}} \put(20,5){\circle{4}} \put(20,75){\circle{4}}
\put(40,40){\circle*{2}} \put(60,40){\circle{4}} \put(60,20){\circle*{2}}
\put(60,60){\circle*{2}} \put(80,40){\circle*{2}} \put(0,40){\line(1,0){18}}
\put(20,42){\line(0,1){31}} \put(20,38){\line(0,-1){31}}
\put(22,40){\line(1,0){36}} \put(60,42){\line(0,1){18}}
\put(60,38){\line(0,-1){18}} \put(62,40){\line(1,0){18}}\end{picture}\] is
rational with Shabat polynomial
$$(x^2-500/441)^4(x^2+500x/189+500/189).$$ Its rotation group has order 1440
and it belong to Galois 1-orbit.

\begin{rem} This tree is, so called, "special" tree, i.e. a tree with a
primitive rotation group with order less, than $n!/2$ ($n$ is the number of
edges) \cite{Adr, LZ}. \ab Remaining 5 trees constitute Galois 5-orbit. Their
rotation group is $S_{10}$.\end{rem}

\subsection{Type 83}  It contains 2 trees:
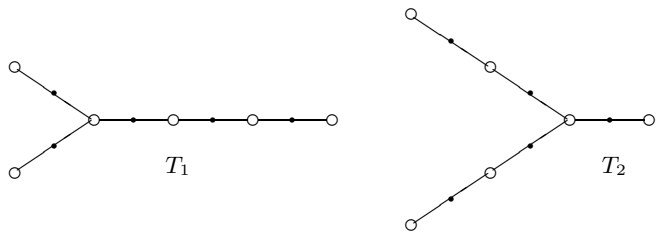
\[\begin{picture}(240,90) \put(0,25){\circle{4}} \put(0,65){\circle{4}}
\put(15,35){\circle*{2}} \put(15,55){\circle*{2}} \put(30,45){\circle{4}}
\put(45,45){\circle*{2}} \put(60,45){\circle{4}} \put(75,45){\circle*{2}}
\put(90,45){\circle{4}} \put(105,45){\circle*{2}} \put(120,45){\circle{4}}
\put(1,26){\line(3,2){28}} \put(1,64){\line(3,-2){28}}
\put(32,45){\line(1,0){26}} \put(62,45){\line(1,0){26}}
\put(92,45){\line(1,0){26}} \put(57,25){\small $T_1$}

\put(150,5){\circle{4}} \put(150,85){\circle{4}} \put(165,15){\circle*{2}}
\put(165,75){\circle*{2}} \put(180,25){\circle{4}} \put(180,65){\circle{4}}
\put(195,35){\circle*{2}} \put(195,55){\circle*{2}} \put(210,45){\circle{4}}
\put(225,45){\circle*{2}} \put(240,45){\circle{4}} \put(151,6){\line(3,2){28}}
\put(151,84){\line(3,-2){28}} \put(181,26){\line(3,2){28}}
\put(181,64){\line(3,-2){28}} \put(212,45){\line(1,0){26}} \put(222,25){\small
$T_2$} \end{picture}\] Both are rational with Shabat polynomials
$$(3x+20)^3(3x-10)^2(x-10)^2(3x^3+20x^2-400x-4000)$$ for tree $T_1$ and
$$(60x-7)^3(1200x^2+280x+343)^2(100x^2+35x+49)(15x-7)$$ for tree $T_2$. Their
rotation groups have orders 14400 and 7200, respectively.

\vspace{1cm}

\begin{thebibliography}{9}
\bibitem{LZ} Lando S.K., Zvonkin A.K. Graphs on surfaces and their
applications, Springer, 2004.
\bibitem{GJ} Goulden I.P., Jackson D.M. The combinatorial relationship between
trees, cacti and certain connection coefficients for the symmetric group,
European J. Combin., 1992, V.13, p. 357-365.
\bibitem{Sc} Schneps L. Dessins d'enfants on the Riemenn sphere, London Math.
Soc. Lecture Notes Series, 1994, V. 200, p. 47-77.
\bibitem{Ko1} Kochetkov Yu. Yu. On non-trivially decomposable types, Russian
Math. Surveys, 1997, V. 52(4), p. 836-837.
\bibitem{Zv} Betrema J., Pere D., Zvonkin A.K. Plane trees and their Shabat
Polynomials. Catalog, Bordeaux: Rapport interne de LaBRI, 1992.
\bibitem{Ko2} Kochetkov Yu.Yu. Plane trees with nine edges. Catalog, J. of
Mathematical Sciences, 2009, V. 158(1), p. 114-140.
\bibitem{Ko3} Kochetkov Yu.Yu. Geometry of plane trees, J. of
Mathematical Sciences, 2009, V. 158(1), p. 106-113.
\bibitem{Adr} Adrianov N.M., Kochetkov Yu.Yu., Suvorov A.D. Plane trees with
exceptional primitive rotation groups [in Russian], Fund. i prikl. mat., 1997,
V. 3(4), p. 1085-1092.
\end{thebibliography}
\end{document}